\documentclass[12pt,leqno,draft]{article}
\usepackage{amsfonts}
\usepackage{amsmath, dsfont, amsthm, amsfonts, amssymb, color}
\usepackage{mathrsfs}
\usepackage{color}
\usepackage{stmaryrd}
\theoremstyle{definition}
\newtheorem{defn}{Definition}[section]
\newcommand{\scr}[1]{\mathscr #1}
\definecolor{wco}{rgb}{0.5,0.2,0.3}

\numberwithin{equation}{section} \theoremstyle{remark}

\newcommand{\ua}{\uparrow}

\title{{\bf 
Derivative Estimates  on Distributions of McKean-Vlasov SDEs }\footnote{Supported in
 part by  NNSFC (11771326, 11831014, 11801406, 11921001).} }
\author{
{\bf   Xing Huang $^{a)}$,  Feng-Yu Wang $^{a), b)}$  }\\
\footnotesize{ a)Center for Applied Mathematics, Tianjin
University, Tianjin 300072, China}\\
\footnotesize{  xinghuang@tju.edu.cn}\\
 \footnotesize{ b)Department of Mathematics,
Swansea University, Singleton Park, SA2 8PP, United Kingdom}\\
\footnotesize{  wangfy@tju.edu.cn}}
\begin{document}
\allowdisplaybreaks
\def\R{\mathbb R}  \def\ff{\frac} \def\ss{\sqrt} \def\B{\mathbf
B} \def\W{\mathbb W}
\def\N{\mathbb N} \def\kk{\kappa} \def\m{{\bf m}}
\def\ee{\varepsilon}\def\ddd{D^*}
\def\dd{\delta} \def\DD{\Delta} \def\vv{\varepsilon} \def\rr{\rho}
\def\<{\langle} \def\>{\rangle} \def\GG{\Gamma} \def\gg{\gamma}
  \def\nn{\nabla} \def\pp{\partial} \def\E{\mathbb E}
\def\d{\text{\rm{d}}} \def\bb{\beta} \def\aa{\alpha} \def\D{\scr D}
  \def\si{\sigma} \def\ess{\text{\rm{ess}}}
\def\beg{\begin} \def\beq{\begin{equation}}  \def\F{\scr F}
\def\Ric{\text{\rm{Ric}}} \def\Hess{\text{\rm{Hess}}}
\def\e{\text{\rm{e}}} \def\ua{\underline a} \def\OO{\Omega}  \def\oo{\omega}
 \def\tt{\tilde} \def\Ric{\text{\rm{Ric}}}
\def\cut{\text{\rm{cut}}} \def\P{\mathbb P} \def\ifn{I_n(f^{\bigotimes n})}
\def\C{\scr C}      \def\aaa{\mathbf{r}}     \def\r{r}
\def\gap{\text{\rm{gap}}} \def\prr{\pi_{{\bf m},\varrho}}  \def\r{\mathbf r}
\def\Z{\mathbb Z} \def\vrr{\varrho} \def\ll{\lambda}
\def\L{\scr L}\def\Tt{\tt} \def\TT{\tt}\def\II{\mathbb I}
\def\i{{\rm in}}\def\Sect{{\rm Sect}}  \def\H{\mathbb H}
\def\M{\scr M}\def\Q{\mathbb Q} \def\texto{\text{o}} \def\LL{\Lambda}
\def\Rank{{\rm Rank}} \def\B{\scr B} \def\i{{\rm i}} \def\HR{\hat{\R}^d}
\def\to{\rightarrow}\def\l{\ell}\def\iint{\int}
\def\EE{\scr E}\def\Cut{{\rm Cut}}
\def\A{\scr A} \def\Lip{{\rm Lip}}
\def\BB{\scr B}\def\Ent{{\rm Ent}}\def\L{\scr L}
\def\R{\mathbb R}  \def\ff{\frac} \def\ss{\sqrt} \def\B{\mathbf
B}
\def\N{\mathbb N} \def\kk{\kappa} \def\m{{\bf m}}
\def\dd{\delta} \def\DD{\Delta} \def\vv{\varepsilon} \def\rr{\rho}
\def\<{\langle} \def\>{\rangle} \def\GG{\Gamma} \def\gg{\gamma}
  \def\nn{\nabla} \def\pp{\partial} \def\E{\mathbb E}
\def\d{\text{\rm{d}}} \def\bb{\beta} \def\aa{\alpha} \def\D{\scr D}
  \def\si{\sigma} \def\ess{\text{\rm{ess}}}
\def\beg{\begin} \def\beq{\begin{equation}}  \def\F{\scr F}
\def\Ric{\text{\rm{Ric}}} \def\Hess{\text{\rm{Hess}}}
\def\e{\text{\rm{e}}} \def\ua{\underline a} \def\OO{\Omega}  \def\oo{\omega}
 \def\tt{\tilde} \def\Ric{\text{\rm{Ric}}}
\def\cut{\text{\rm{cut}}} \def\P{\mathbb P} \def\ifn{I_n(f^{\bigotimes n})}
\def\C{\scr C}      \def\aaa{\mathbf{r}}     \def\r{r}
\def\gap{\text{\rm{gap}}} \def\prr{\pi_{{\bf m},\varrho}}  \def\r{\mathbf r}
\def\Z{\mathbb Z} \def\vrr{\varrho} \def\ll{\lambda}
\def\L{\scr L}\def\Tt{\tt} \def\TT{\tt}\def\II{\mathbb I}
\def\i{{\rm in}}\def\Sect{{\rm Sect}}  \def\H{\mathbb H}
\def\M{\scr M}\def\Q{\mathbb Q} \def\texto{\text{o}} \def\LL{\Lambda}
\def\Rank{{\rm Rank}} \def\B{\scr B} \def\i{{\rm i}} \def\HR{\hat{\R}^d}
\def\to{\rightarrow}\def\l{\ell}
\def\8{\infty}\def\I{1}\def\U{\scr U}
\maketitle

\begin{abstract} By using the heat kernel parameter expansion   with respect to the frozen SDEs, the intrinsic derivative is estimated for the law of
 Mckean-Vlasov SDEs with respect to the initial distribution.  As an application,
 the total variation distance between the laws of two solutions  is bounded by the Wasserstein distance for initial distributions. These extend some recent  results
 proved  for  distribution-free noise     by using  the coupling  method  and Malliavin calculus.
 \end{abstract} \noindent
 AMS subject Classification:\  60H1075, 60G44.   \\
\noindent
 Keywords: Mckean-Vlasov SDEs, intrinsic derivative, $L$-derivative,  heat kernel parameter expansion.
 \vskip 2cm

\section{Introduction}
Let $\scr P_2$ be the set of all probability measures on $\R^d$ with finite second moment, which is called the Wasserstein space under the metric
$$\W_2(\mu,\nu):= \inf_{\pi\in \C(\mu,\nu)} \bigg(\int_{\R^d\times\R^d} |x-y|^2 \pi(\d x,\d y)\bigg)^{\ff 1 {2}},\ \ \mu,\nu\in \scr P_{2},$$ where $\C(\mu,\nu)$ is the set of all couplings of $\mu$ and $\nu$.   Consider the following distribution dependent SDE on $\R^d$:
\beq\label{E1} \d X_t^\mu= b_t(X_t^\mu , \L_{X_t^\mu})\d t+\si_t(X_t^\mu,\L_{X_t^\mu})\d W_{t},\ \ \L_{X_0^\mu}=\mu\in\scr P_2,\end{equation}
where $W_t$ is an $m$-dimensional Brownian motion on a complete filtration probability space $(\OO,\{\F_t\}_{t\ge 0},\P)$, $\L_{X_t}$ is the law of $X_t$, and
\begin{align*} &b: \R_+\times\R^d\times \scr P_2\to \R^d,\ \ \si: \R_+\times \R^d\times \scr P_2\to \R^d\otimes\R^m
 \end{align*}
 are measurable.  This type equations, known as  Mckean-Vlasov or mean field SDEs, have been intensively investigated and applied, see for instance the monograph \cite{CDLL19}
 and references therein.

 To characterize  the regularity of   the law $\L_{X_t^\mu}$ with respect to the initial distribution $\mu$,  we  investigate  the derivative estimate of the functions
 $$\scr P_2\ni\mu\mapsto P_tf(\mu):= \E f(X_t^\mu),\ \ f\in \B_b(\R^d), t>0.$$
 When  the noise coefficient $\si_t(x,\mu)$ does not depend on $\mu$, the Harnack inequality and derivative formula have been established in \cite{FYW1, RW}
 for   $P_tf $    by using the coupling by change of measures and Malliavin calculus respectively. See also \cite{BRW20, HRW, HW, Song} for extensions to distribution-path dependent SDEs/SPDEs,  singular distribution dependent SDEs, and  distribution dependent SDEs with jumps, where in \cite{Song} allows the noise to be also  distribution dependent and establishes the gradient estimate on $P_t f(x):=(P_tf)(\dd_x)$ when the initial distribution is a Dirac measure.  In this paper, we   estimate the derivative
 of $P_tf(\mu)$  in $\mu$ by using the heat kernel parameter expansion  with respect to the frozen SDE
 \beq\label{E1'} \d X_t^{z,\mu} = b_t(z,  \mu_t)\d t+\si_t(z,\mu_t)\d W_{t}\end{equation}
 for fixed $(z,\mu)\in\R^d\times \scr P_2(\R^d),$ where  $\mu_t:= \L_{X_t^\mu}.$  Since this SDE has constant coefficients, the solution has a Gaussian heat kernel which can be easily analyzed.

Before introducing the main result, we first recall the intrinsic derivative and $L$-derivative for functions on $\scr P_2$ which go back to \cite{AKR} where the intrinsic derivative on the configurations space is introduced, see \cite{RW19} for the link of different derivatives for measures.

\begin{defn}\label{defLio} Let $f:\scr P_2\to \mathbb{R}$ and $g:\mathbb{R}^d\times \scr P_2\to \mathbb{R}$.
\begin{enumerate}
\item[(1)] $f$ is called intrinsically differentiable, if for any $\mu\in \scr P_2$,
$$L^2(\R^d\to\R^d;\mu)\ni\phi\mapsto D_\phi^Lf(\mu):=\lim_{\vv\downarrow 0} \ff{f(\mu\circ({\rm Id}+\vv\phi)^{-1})- f(\mu)}\vv\in\R$$
is a well defined bounded linear functional. In this case,  the    unique map
$$\scr P_2\ni\mu\mapsto D^Lf(\mu)\in L^2(\R^d\to\R^d;\mu)$$
such that $D_\phi^Lf(\mu)= \<\phi, D^Lf(\mu)\>_{L^2(\mu)}$ holds for any $\mu\in \scr P_2$  and $\phi\in L^2(\R^d\to\R^d;\mu)$ is called the intrinsic derivative of $f$, and we denote
$$\|D^Lf(\mu)\|:= \|D^Lf(\mu)(\cdot)\|_{L^2(\mu)},\ \ \mu\in \scr P_2.$$
If moreover
$$\lim_{\mu(|\phi|^2)\to 0}\frac{f(\mu\circ({\rm Id}+ \phi)^{-1})-f(\mu)-D_\phi^Lf(\mu)}{\sqrt{\mu(|\phi|^2)}}=0,\ \ \mu\in \scr P_2,$$ we call $f$ $L$-differentiable, and in this case $D^Lf$ is also called the $L$-derivative of $f$.
\item[(2)] We denote $f\in C^1(\scr P_2)$, if $f$ is $L$-differentiable and  its $L$-derivative has a version $D^L f(\mu)(x)$ jointly continuous in $(x,\mu)\in \R^d\times \scr P_2.$
 \item[(3)] We denote $g\in C^{1,1}(\R^d\times\scr P_2)$, if $g(x,\cdot)\in C^1(\scr P_2)$ for $x\in \R^d$, $g(\cdot,\mu)\in C^1(\R^d)$ for $\mu\in \scr P_2$,
 $g(x,\mu), \nn g(\cdot,\mu)(x)$ are jointly continuous in $(x,\mu)\in \R^d\times\scr P_2$,  and $D^Lg(x,\cdot)(\mu)(y)$  has a version jointly continuous in $(x,y,\mu)\in \R^d\times\R^d\times\scr P_2.$
 \item[(4)] A vector- or matrix-valued function is said in a class defined above, if so are its component functions.
 \end{enumerate}
\end{defn}

To estimate the intrinsic derivative of $P_tf(\mu)$, we need the following condition. Let $|\cdot|$ and $\|\cdot\|$ denote the norm in $\R^d$ and the operator norm for linear operators repsectively.

 \beg{enumerate} \item[{\bf (H)}] For any $t\ge 0$, $b_t,\si_t\in C^{1,1}(\R^d\times\scr P_2),$  and  there exists an increasing function $K: [0,\infty)\to [0,\infty)$ such that
 for any  $ t\ge 0, x,y\in\R^d$ and $ \mu\in \scr P_2(\R^d),$
$$  K_t^{-1}  {\rm Id} \le  (\si_t\si_t^*) (x,\mu)\le  K_t {\rm Id},$$
\beg{align*}  &\ |b_t(x,\mu)|+ \|\nn b_t(\cdot,\mu)(x)\|+ \|D^L \{b_t(x,\cdot)\}(\mu)\| \\
&\qquad  +  \|\nn \{\si_t(\cdot,\mu)\}(x)\|^2+  \|D^L\{\si_t(x,\cdot)\}(\mu)\|^2
  \le K_t,\end{align*}
\beg{align*}   &\|D^L \{b_t(x,\cdot)\}(\mu)- D^L \{b_t(y,\cdot)\}(\mu)\|+ \|D^L \{\si_t(x,\cdot)\}(\mu)- D^L\{ \si_t(y,\cdot)\}(\mu)\|\\
&\le K_t |x-y|^2.
  \end{align*}
\end{enumerate}

 It is well known that    SDE \eqref{E1} is well-posed under the assumption {\bf (H)}, so that $P_tf$ is  well defined on $\scr P_2$ for any $t\ge 0$ and $f\in \B_b(\R^d).$
 In general, for any $s\ge 0$ and $X_{s,s}^\mu\in L^2(\OO\to\R^d, \F_s,\P)$ with $\L_{X_{s,s}^\mu} =\mu$, let $X_{s,t}^\mu$ be the unique solution of \eqref{E1} for $t\ge s$:
 \beq\label{E11} \d X_{s,t}^\mu= b_t(X_{s,t}^\mu , \L_{X_{s,t}^\mu})\d t+\si_t(X_{s,t}^\mu,\L_{X_{s,t}^\mu})\d W_{t},\ \ t\ge s, \L_{X_{s,s}^\mu}=\mu\in\scr P_2.\end{equation}
 We denote  $P_{s,t}^*\mu=\L_{X_{s,t}^\mu}$ and investigate the regularity of
 $$P_{s,t}f(\mu):= \E f(X_{s,t}^\mu)= \int_{\R^d} f\d (P_{s,t}^*\mu),\ \ f\in \B_b(\R^d).$$
 By the uniqueness, we have the flow property
 $$P_{s,t}^*= P_{r,t}^*P_{s,r}^*,\ \ 0\le s\le r\le t.$$ However,  due to the distribution dependence,  $P_{s,t}$ is no-longer a semigroup, i.e. in general
 $P_{s,t}\ne P_{r,t}P_{s,r}$  and
 $$P_tf(\mu)\ne \int_{\R^d} P_tf(x)\mu(\d x),$$ so that the regularity of $P_tf(\mu)$ in $\mu\in\scr P_2$ can not be  deduced from that of $P_tf(x):= P_tf(\dd_x)$ for $x\in \R^d$, see for instance \cite{FYW1} for details.

We now state  the main result of the paper as follows.

 \beg{thm}\label{T1} Assume {\bf (H)}.  Then for any $t>s$ and $f\in \B_b(\R^d)$, $P_{s,t}f$ is
  $L$-differentiable, and there exists an increasing function $C: [0,\infty)\to (0,\infty)$ such that
\beq\label{EST1} \|D^L P_{s,t} f(\mu)\|\le \ff{C_t \|f\|_\infty}{\ss {t-s}},\ \ t>s\ge 0, f\in \B_b(\R^d).\end{equation}
 Consequently, for any $t>s\ge 0,   \mu,\nu\in \scr P_2,$
\beq\label{EST2} \|P_{s,t}^*\mu- P_{s,t}^*\nu\|_{var}:= \sup_{\|f\|_\infty\le 1} |P_{s,t}f(\mu)-P_{s,t}f(\nu)|\le \ff{C_t\|f\|_\infty}{\ss {t-s}} \W_2(\mu,\nu).\end{equation}
\end{thm}

\paragraph{Remark 1.1.} We may also apply Malliavin calculus to establish a derivative formula for $D^L P_{s,t} f(\mu)$  as in \cite{Song}, where the usual derivative in initial points
(rather than in initial distributions) are studied. However, in this way we need stronger conditions on the coefficients,  i.e.    $b_t(x,\mu)$ and $\si_t(x,\mu)$ also have bounded second order derivatives in $x$. Let us explain this in more details.

Firstly, under  {\bf (H)},  the Malliavin matrix
  $$M_{s,t}:= \big\{\<D (X_{s,t}^\mu)_i, D (X_{s,t}^\mu)_j\>_\H\big\}_{1\le i\le j}$$
  is invertible with $\E\|M_{s,t}^{-1}\|^2<\infty$ for $t>s\ge 0$, where $D$ is the Malliavin gradient,   $\H$ is the Cameron-Martin space in Malliavin calculus, and $(X_{s,t}^\mu)_i$ is   the $i$-th component of $X_{s,t}^\mu$.

  Next, for any $\phi\in L^2(\R^d\to\R^d;\mu)$, let $v_{s,t}^\phi= D_\phi^L X_{s,t}^\mu$, which exists  in $L^2(\P)$ and satisfies
  $$\E|v_{s,t}^\phi|^2\le c(t) \mu(|\phi|^2)$$ for some constant $c(t)>0$, see \cite[Proposition 3.2]{RW}.

  Then for any $f\in C_b^1(\R^d)$, by the chain rule and the integration by parts formula for
  the Malliavin gradient $D$, we have
  \beg{align*} & D_\phi^LP_{s,t}f(\mu)= \E\<\nn f(X_{s,t}^\mu), v_{s,t}^\phi\> = \sum_{i=1}^d \E[\pp_i f(X_{s,t}^\mu) (v_{s,t}^\phi)_i] \\
  &= \sum_{i,j,k=1}^d\E\big[ \pp_i f(X_{s,t}^\mu) (M_{s,t})_{ij} (M_{s,t}^{-1})_{jk} (v_{s,t}^\phi)_k\big]\\
  &=  \sum_{i,j,k=1}^d\E\big[\<D f(X_{s,t}^\mu), D(X_{s,t}^\mu)_j\>_\H    (M_{s,t}^{-1})_{jk} (v_{s,t}^\phi)_k\big]\\
  &= \sum_{i,j,k=1}^d \E\big[f(X_{s,t}^\mu) D^*\big\{ (M_{s,t}^{-1})_{jk} (v_{s,t}^\phi)_k D(X_{s,t}^\mu)_j\big\}\big],\end{align*}
  where $D^*$ is the Malliavin divergence. To make the above calculations meaningful, we need to verify that $(M_{s,t}^{-1})_{jk} (v_{s,t}^\phi)_k D(X_{s,t}^\mu)_j$ belongs to the domain of   $D^*$, for which the second order derivatives of coefficients will be involved. For instance, as shown in \cite[Proposition 3.2]{RW} that $v_{s,t}^\phi$ solves an SDE
  involving in the first order derivatives of $b$ and $\si$,  making Malliavin derivative to this SDE we see that $D v_{s,t}^\phi$  solves an SDE containing  the second order derivatives of coefficients.

  \

The remainder of the paper is organized as follows. In Section 2,  we formulate $P_{s,t}f(\mu)$ using classical SDEs with parameter $\mu$ and the parameter expansion of heat kernels with respect to the  frozen SDE \eqref{E1'},  and estimate the $L$-derivative for functions of $P_{s,t}^*\mu$.   With these preparations, we prove Theorem \ref{T1} in Section 3.

 \section{Preparations}

 We first represent $P_{s,t}f(\mu)$ by using a Markov semigroup $P_{s,t}^\mu$ with parameter $\mu$, then introduce the heat kernel expansion of $P_{s,t}^\mu$ with respect to the frozen SDEs. Since the frozen SDE has explicit Gaussian heat kernel, this enables us to calculate the intrinsic derivative of $P_tf(\mu)$ with respect to $\mu$.

\subsection{A representation of $P_{s,t}$}

   For any $s\ge 0, x\in\R^d$ and $\mu\in \scr P_2$, consider the decoupled SDE
   \beq\label{ES0} \d X_{s,t}^{x,\mu}= b_t(X_{s,t}^{x,\mu}, P_{s,t}^*\mu)\d t + \sigma_t(X_{s,t}^{x,\mu}, P_{s,t}^*\mu)\d W_t, \ \ X_{s,s}^{x,\mu}=x, t\ge s.\end{equation}
   In this SDE,  the measure variable $P_{s,t}^*\mu$ is fixed, so that it reduces to the classical time inhomogeneous SDE.  Let $P_{s,t}^\mu$ be the associated Markov semigroup, i.e.
   $$P_{s,t}^\mu f(x)= \E f(X_{s,t}^{x,\mu}),\ \ t\ge s, f\in \B_b(\R^d), x\in \R^d.$$
Since  $X_{s,t}^\mu$ solves \eqref{ES0} with the random initial value    $X_{s,s}^\mu$ replacing $x$,  and since $\L_{X_{s,s}^\mu}=\mu,$  by the standard Markov property of solutions to \eqref{ES0}, we have
\beq\label{ES1} P_{s,t}f(\mu):=\E f(X_{s,t}^\mu)= \int_{\R^d} P_{s,t}^\mu f(x) \mu(\d x),\ \ t\ge s, f\in \B_b(\R^d), \mu\in \scr P_2.\end{equation}
Since for any $g\in C_b^1(\R^d)$ the function
$\mu\mapsto \mu(g):= \int_{\R^d}g\d\mu$ is $L$-differentiable with
$D^L \mu(g)= \nn g$,  we first study the derivative of $P_{s,t}^\mu f(x)$ in $x$.

\beg{lem}\label{L0} Assume {\bf (H)}. Then for any $f\in \B_b(\R^d)$ and $t>s\ge 0$,
we have $P_{s,t}^\mu f\in C^1(\R^d)$ such that $(\nn P_{s,t}^\mu f)(x)$  is continuous in $(x,\mu)\in \R^d\times\scr P_2$, and
\beq\label{ES2} \|\nn P_{s,t}^\mu f\|_\infty \le \ff{CK_t \|f\|_\infty }{\ss{t-s}} \e^{C K_t},\ \ t>s, f\in \B_b(\R^d), \mu\in \scr P_2\end{equation}
holds for some constant $C>0.$\end{lem}

\beg{proof} Since {\bf (H)} implies that $P_{s,t}^*\mu$ is Lipschitz continuous in $\mu\in \scr P_2$, see for instance \cite{FYW1}, the desired assertions follow from    {\bf (H)} and the   Bismut formula
\beq\label{BS} \nn_v P_{s,t}^\mu f(x)  =  \E\bigg [\ff{f(X_{s,t}^{x,\mu})}{t-s}\int_s^t \big\<\{\si_r (\si_r\si_r^*)^{-1}\}(X_{s,r}^{x,\mu}, P_{s,r}^*\mu)  v_{s,r}^{x,\mu}, \d W_r\big\>\bigg],\ \ v\in\R^d\end{equation} for $ f\in \B_b(\R^d)$,
where $v_{s,t}^{x,\mu}:=\ff{\d}{\d \vv} X_{s,t}^{x+\vv v, \mu}|_{\vv=0} $ solves the linear SDE
\beq\label{LL}\d v_{s,t}^{x,\mu}= \{\nn_{v_{s,t}^{x,\mu}}b_t(\cdot, P_{s,t}^*\mu) \}(X_{s,t}^{x,\mu} )\d t +\{\nn_{v_{s,t}^{x,\mu}} \si_t(\cdot,P_{s,t}^*\mu)\}(X_{s,t}^{x,\mu}) \d W_t,
\ \ t\ge s, v_{s,s}^{x,\mu}=v.\end{equation}  By {\bf (H)},   $v_{s,t}^{x,\mu}$ is continuous in $(x,\mu)\in \R^d\times\scr P_2$ and
$$\E |v_{s,t}^{x,\mu}|^2\le |v|^2 \e^{C K_t},\ \ t\ge s, v\in \R^d$$ holds for some constant $C>0$, so that \eqref{BS} implies that $(\nn P_{s,t}^\mu f)(x)$  is continuous in $(x,\mu)\in \R^d\times\scr P_2$  and satisfies \eqref{ES2}.

To prove \eqref{BS}, for fixed $t>s$, take
$$h_u=\int_s^u \{\si_r^* (\si_r\si_r^*)^{-1}\}(X_{s,r}^{x,\mu}, P_{s,r}^*\mu)  v_{s,r}^{x,\mu} \d r,\ \ u\in [s,t].$$   Then the Malliavin derivative $w_r^{x,\mu}:= D_h X_{s,r}^{x,\mu}$ along $h$ solves the SDE
\beg{align*} &\d w_{s,r}^{x,\mu} = \Big[\{\nn_{w_{s,r}^{x,\mu}}b_r(\cdot, P_{s,r}^*\mu) \}(X_{s,r}^{x,\mu} )+\si_r(X_{s,r}^{x,\mu}, P_{s,r}^*\mu) h_r'\Big] \d r +\{\nn_{w_{s,r}^{x,\mu}} \si_r(\cdot,P_{s,r}^*\mu)\}(X_{s,r}^{x,\mu}) \d W_r\\
&=\Big[\{\nn_{w_{s,r}^{x,\mu}}b_r(\cdot, P_{s,r}^*\mu) \}(X_{s,r}^{x,\mu} )+ v_{s,r}^{x,\mu} \Big] \d r +\{\nn_{w_{s,r}^{x,\mu}} \si_r(\cdot,P_{s,r}^*\mu)\}(X_{s,r}^{x,\mu}) \d W_r ,\ \ r\in [s,t],  w_{s,s}^{x,\mu}=0,\end{align*} see for instance \cite[Proposition 3.5]{RW}.
It is easy to see from \eqref{LL} that $\bar v_r:= (r-s)v_{s,r}^{x,\mu}$ solves the same equation. By the uniqueness we obtain
$(t-s)v_{s,t}^{x,\mu}= D_h X_{s,t}^{x,\mu}$, so that the chain rule and the integration by parts formula yield
\beg{align*} &\nn_v P_{s,t}^\mu f(x)= \E\<\nn f(X_{s,t}^{x,\mu}), v_{s,t}^{x,\mu}\>=\ff 1 {t-s} \E\<\nn f(X_{s,t}^{x,\mu}), D_h X_{s,t}^{x,\mu}\>\\
&= \ff 1 {t-s} \E D_h\{ f(X_{s,t}^{x,\mu})\} =   \E\bigg[ \ff{f(X_{s,t}^{x,\mu})}{t-s} \int_s^t \big\<\{\si_r^\ast (\si_r\si_r^*)^{-1}\}(X_{s,r}^{x,\mu}, P_{s,r}^*\mu)  v_{s,r}^{x,\mu}, \d W_r\big\>\bigg].\end{align*}
\end{proof}

Combining \eqref{ES1} with Lemma \ref{L0}, we have the following result.

\beg{lem}\label{L1} Assume {\bf (H)}. Let $t>s$ and $f\in \B_b(\R^d)$. If for any $x\in\R^d$, the function $\mu\mapsto P_{s,t}^\mu f(x)$ is $L$-differentiable with
\beq\label{ES3}  \sup_{x\in \R^d} \big\|D^L \{P_{s,t}^\cdot f(x)\}(\mu)\big\|<\infty,\end{equation}
then $P_{s,t}f(\mu)$ is $L$-differentiable in $\mu$ with
\beq\label{ES4} D^L P_{s,t}f(\mu)= \nn P_{s,t}^\mu f + \int_{\R^d} D^L \{P_{s,t}^\cdot f(x)\}(\mu) \mu(\d x).\end{equation}
Consequently,  there exists a constant $C>0$ such that for any $f\in \B_b(\R^d)$ and $\mu\in\scr P_2,$
\beq\label{ES5} \big\|D^L P_{s,t}f(\mu)\big\|\le \ff{CK_t \|f\|_\infty }{\ss{t-s}} \e^{C K_t}+ \sup_{x\in \R^d} \big\|D^L \{P_{s,t}^\cdot f(x)\}(\mu)\big\|, \ \ t>s\ge 0.\end{equation}
\end{lem}

\beg{proof} Obviously, \eqref{ES5} is implied by \eqref{ES2} and \eqref{ES4}. So, we only need to prove that $P_{s,t}f(\mu)$ is $L$-differentiable and satisfies \eqref{ES4}.

(1) We first prove that  $P_{s,t}f(\mu)$ is  intrinsically differentiable and satisfies \eqref{ES4}.   For any $g\in C_b^1(\R^d)$, the function
$\mu\mapsto \mu(g):= \int_{\R^d}g\d\mu$ is $L$-differentiable with
$D^L \mu(g)= \nn g$. So, for any $\mu\in \scr P_2$, the function
\begin{align}\label{Cy}\scr P_2\ni\nu\mapsto  P_{s,t}^\mu f(\nu):= \int_{\R^d} P_{s,t}^\mu f\d\nu
 \end{align}
 is $L$-differentiable with
$D^L (P_{s,t}^\mu f)(\nu)= \nn P_{s,t}^\mu f,\ \nu\in\scr P_2.$ Combining this with \eqref{ES1}, \eqref{ES2} and \eqref{ES3},  and using the dominated convergence theorem, we conclude that the map
$$L^2(\R^d\to\R^d;\mu)\ni \phi\mapsto  D^L_\phi P_{s,t}f(\mu)= \<\nn P_{s,t}^\mu f,\phi\>_{L^2(\mu)} + \int_{\R^d} D^L_\phi \{P_{s,t}^\cdot f(x)\}(\mu) \mu(\d x)$$
is a bounded linear functional, so that by definition,
$P_{s,t}f(\mu)$ is  intrinsically differentiable in $\mu\in \scr P_2$, and the formula \eqref{ES4} holds true.

(2) By \eqref{ES4}, for any $\phi\in L^2(\R^d\to\R^d;\mu)$, we have
\beg{align*}& P_{s,t}f (\mu\circ({\rm Id}+\phi)^{-1} )- P_{s,t}f (\mu )- D_\phi^L P_{s,t}f(\mu)\\
&= \int_{\R^d} \big\{P_{s,t}^{\mu\circ ({\rm Id}+\phi)^{-1}} f(x+\phi(x)) - P_{s,t}^{\mu\circ ({\rm Id}+\phi)^{-1}} f(x) - \<\nn P_{s,t}^{\mu\circ ({\rm Id}+\phi)^{-1}} f(x), \phi(x)\>\big\}\mu(\d x)\\
& \quad + \int_{\R^d} \big\{P_{s,t}^{\mu\circ ({\rm Id}+\phi)^{-1}} f(x)- P_{s,t}^\mu f(x) - D_\phi^L[P_{s,t}^\cdot f(x)](\mu) \big\}\mu(\d x)\\
&\quad +\int_{\R^d} \<\nn P_{s,t}^{\mu\circ ({\rm Id}+\phi)^{-1}} f(x)- \nn P_{s,t}^\mu f(x), \phi(x)\>\mu(\d x).\end{align*}
Combining this with Lemma \ref{L0}, \eqref{ES3}, and the $L$-differentiability of $P_{s,t}^\mu f(x)$  in $\mu$, we may apply the dominated convergence theorem to derive
$$\lim_{\|\phi\|_{L^2(\mu)}\downarrow 0} \ff{|P_{s,t}f (\mu\circ({\rm Id}+\phi)^{-1} )- P_{s,t}f (\mu )-D_\phi^L P_{s,t}f(\mu)|}{\|\phi\|_{L^2(\mu)}} =0,$$
that is, $P_{s,t}f(\mu)$ is $L$-differentiable.
\end{proof}

 According to Lemma \ref{L1},  to estimate $\|D^L P_{s,t} f(\mu)\|,$ it remains to investigate the  $L$-derivative of $P_{s,t}^\mu f(x)$ in $\mu$. To this end, we let $p_{s,t}^\mu(x,y)$ be the  heat kernel of $P_{s,t}^\mu$ for $t>s$, which exists and is differentiable in $x$ and $y$  under conditions {\bf (H)}. We have
\beq\label{MA0} P_{s,t}^\mu f(x)= \int_{\R^d} p_{s,t}^\mu (x,y) f(y) \d y,\ \ f\in \B_b(\R^d), t>s, x\in \R^d.\end{equation}
 So, to investigate the   $L$-derivative of $P_{s,t}^\mu f(x),$ we need to study  that of $p_{s,t}^\mu(x,y)$, for which we will use the heat kernel parameter expansion.

 \subsection{Parameter expansion for $p_{s,t}^\mu$}

Since  heat kernel  $p_{s,t}^\mu$ is less explicit, we make use of its parameter expansion with respect to the heat kernel of the Gaussian process
 $$X_{s,r,t}^{x,\mu,z}= x + \int_r^t b_u(z, P_{s,u}^*\mu)\d u + \int_r^t \si_u(z, P_{s,u}^*\mu)\d W_u,\ \ t\ge r\ge s\ge 0, x\in\R^d$$
 for fixed $z\in \R^d$ and $\mu\in \scr P_2$.  For any $t\ge r\ge s\ge 0$, let
\beq\label{MA}\beg{split} & m_{s,r,t}^{\mu,z} := \int_r^t b_u(z, P_{s,u}^*\mu)\d u,\ \ m^{\mu,z}_{s,t}:=m^{\mu,z}_{s,s,t},\\
& \ a_{s,r,t}^{\mu,z} :=  \int_r^t (\si_u\si_u^*)(z, P_{s,u}^*\mu)\d u, \ \ a^{\mu,z}_{s,t}:=a^{\mu,z}_{s,s,t}.\end{split}\end{equation}
By {\bf (H)}, we have
\begin{align}\label{mv}
|m_{s,r,t}^{\mu,z}|+|a_{s,r,t}^{\mu,z}|\leq (t-r)K_t,\ \ t\ge r\ge s\ge 0.
\end{align}
Obviously,  the law of $X_{s,r,t}^{x,\mu,z}$ is the $d$-dimensional normal distribution entered at $x+m_{s,r,t}^{\mu,z}$ with covariance matrix $a_{s,r,t}^{\mu,z},$ i.e.  the distribution density function   is
\beq\label{ES6} p_{s,r,t}^{\mu,z}(x,y)= \ff{\exp[-\ff 1 2 \<(a_{s,r,t}^{\mu,z})^{-1}(y-x-m_{s,r,t}^{\mu,z}),  y-x-m_{s,r,t}^{\mu,z}\>]}{(2\pi)^{\ff d 2} ({\rm det} \{a_{s,r,t}^{\mu,z}\})^{\ff 1 2}},\ \ y\in \R^d, t>r\ge s.\end{equation}
When $r=s$, we simply denote $p_{s,t}^{\mu,z}=p_{s,s,t}^{\mu,z}$, so that
\beq\label{ES6'} p_{s,t}^{\mu,z}(x,y)= \ff{\exp[-\ff 1 2 \<(a_{s,t}^{\mu,z})^{-1}(y-x-m_{s,t}^{\mu,z}),  y-x-m_{s,t}^{\mu,z}\>]}{(2\pi)^{\ff d 2} ({\rm det} \{a_{s,t}^{\mu,z}\})^{\ff 1 2}},\ \ y\in \R^d, t>  s.\end{equation}
For any $0\le s\le r<t$ and $y,z\in\R^d$, let
\beq\label{ES7}\beg{split}  H_{s,r,t}^{\mu} (y,z):=  & \left\<b_r (z, P_{s,r}^*\mu) -b_{r}(y, P_{s,r}^*\mu),\nabla p_{s,r,t}^{\mu,z}(\cdot,z)(y)\right\>\\
&+\frac{1}{2}\mathrm{tr}\left[\left\{(\sigma_{r}\sigma^\ast_{r})(z, P_{s,r}^*\mu)-(\sigma_{r}\sigma^\ast_{r})(y, P_{s,r}^*\mu)\right\}\nabla^2 p_{s,r,t}^{\mu,z} (\cdot,z)(y)\right].
\end{split} \end{equation}
By the parameter expansion, see for instance \cite[Lemma 3.1]{KM00},  we have
\beq\label{ES8} p_{s,t}^\mu(x,z) = p_{s,t}^{\mu,z}(x,z) + \sum_{m=1}^\infty \int_s^t \d r \int_{\R^d} H_{s,r,t}^{\mu,m}  (y,z) p_{s,r}^{\mu,z} (x,y) \d y,\end{equation}
where   $H_{s,r,t}^{\mu,m}$ for $m\in \mathbb N$ are defined by
\beq\label{HU} \beg{split}& H_{s,r,t}^{\mu, 1}:= H_{s,r,t}^{\mu},\\
&H_{s,r,t}^{\mu, m} (y,z) := \int_r^t\d u \int_{\R^d}  H_{s,u,t}^{\mu, m-1}(z',z) H_{s,r,u}^{ \mu}(y,z') \d z',\ m\ge 2.\end{split}\end{equation}
Combining  \eqref{ES8}  with \eqref{MA}, \eqref{ES6} and \eqref{ES6'},  to estimate $D^L P_{s,t}^\mu f$,  it suffices to study the $L$-derivative of  $b_r(y, P_{u_1,u_2}^* \mu)$ and
$ (\si_r\si_r^*)(y, P_{u_1,u_2}^* \mu)$ in $\mu$ for $r\ge 0$ and $u_2\ge u_1\ge 0$.
So, we present the following lemma.

\beg{lem}\label{L2}  Assume {\bf (H)} and let $t>s\ge 0$. Then for any $F\in C^1(\scr P_2)$ with bounded $\|D^L F\|$,  $F(P_{s,t}^*\mu)$ is $L$-differentiable in $\mu$ such that
\beq\label{LL1} \|D^L F(P_{s,t}^*\cdot)(\mu)\|\le \|D^L F\|_\infty \e^{4 K_t (t-s)}.\end{equation} Consequently,
for any $r\ge 0, t\ge s\ge 0$ and $y\in \R^d$, $b_r(y, P_{s,t}^* \mu)$ and $(\si_r\si_r^*)(y,P_{s,t}^*\mu)$ are   $L$-differentiable in $\mu$, and
$$\max\Big\{\|D^L b_r(y, P_{s,t}^*\cdot)(\mu)\|,  \|D^L (\si_r\si_r^*)(y, P_{s,t}^*\cdot)(\mu)\|\Big\}\le K_r\e^{4K_t(t-s)}, \ \  \mu\in \scr P_2.$$
\end{lem}

\beg{proof} It  suffices  to prove the first assertion.
We first prove the intrinsic differentiability.  Let $\mu\in \scr P_2$ and $\phi\in L^2(\R^d\to\R^d;\mu).$  Since  $\L_{X_{s,s}^\mu}=\mu$ implies
$$\L_{X_{s,s}^\mu+\vv \phi(X_{s,s}^\mu)}=\mu\circ({\rm Id}+\vv\phi)^{-1},\ \  \vv\ge 0,$$ we have $\L_{X_{s,t}^\vv}= P_{s,t}^*(\mu\circ({\rm Id}+\vv\phi)^{-1})$ for $X_{s,t}^\vv$ solving \eqref{E11}   with initial value  $X_{s,s}^\vv= X_{s,s}^\mu+\vv \phi(X_{s,s}^\mu).$
By \cite[Proposition 3.1]{RW} for $\eta=\phi(X_0^\mu)$ and \cite[(4.21)]{RW} for time $s$ replacing $0$,  for any $\dd\ge 0$,
$$v_{s,t}^{\phi,\dd}:= D_\phi^L X_{s,t}^\dd= \lim_{\vv\downarrow 0} \ff{X_{s,t}^{\dd+\vv}-X_{s,t}^\dd}\vv,\ \ t\ge s$$
exists in  $L^2(\OO\to C([s,T];\R^d);\P)$ for any $T>0$, and solves the linear SDEs:
\begin{equation}\begin{split}\label{vp}
\d v_{s,t}^{\phi,\dd}&= \Big[\nabla_{v_{s,t}^{\phi,\dd}}b_t(X^\dd_{s,t}, \L_{X^\dd_{s,t}})+\E \big\{\<D^Lb_t(z,\cdot)( \L_{X^\dd_{s,t}})(X^\dd_{s,t}), v_{s,t}^{\phi,\delta}\>\big\}\big|_{z=X^\dd_{s,t}}\Big]\d t\\
&\quad +\Big[\nabla_{v_{s,t}^{\phi,\dd}}\sigma_t(X^\dd_{s,t}, \L_{X^\dd_{s,t}})+\E \big\{D^L\sigma_t(z,\cdot)( \L_{X^\dd_{s,t}})(X^\dd_{s,t}) v_{s,t}^{\phi,\dd}\big\}\big|_{z=X^\dd_{s,t}}\Big]\d W_{t},\\
&\quad    v_{s,s}^{\phi,\dd}=\phi(X_0),\ \ t\ge s.
\end{split}\end{equation}
Fromm  {\bf (H)} we see that $v_{s,t}^{\phi,\vv}$ is continuous in $\vv$ and
\beq\label{ES9} \E |v_{s,t}^{\phi,\dd}|^2\le \mu(|\phi|^2) \e^{8 (t-s) K_t},\ \ t\ge s, \phi\in L^2(\R^d\times\R^d;\mu).\end{equation}
By the chain rule, see for instance \cite[Proposition 3.1]{RW}, we have
\beq\label{LL2} D^L_\phi F( P_{s,t}^*\cdot)(\mu) =\ff{\d}{\d\vv} F(\L_{X_{s,t}^\vv})\Big|_{\vv=0}= \E\<(D^L F)(P_{s,t}^*\mu)(X_{s,t}^\mu), v_{s,t}^{\phi,0}\>.\end{equation}
Combining this with {\bf(H)} and \eqref{ES9},  we obtain
\beg{align*} &|D^L_\phi  F( P_{s,t}^*\cdot)(\mu)|\le \|(D^L F)(P_{s,t}^*\mu)\|\ss{ \E |v_{s,t}^{\phi,0}|^2}\\
&\le \|\phi\|_{L^2(\mu)}  \|D^LF\|_\infty  \e^{4(t-s) K_t},\ \ \phi\in L^2(\R^d\to\R^d;\mu).\end{align*}
Therefore, $ F(P_{s,t}^* \mu)$ is intrinsically differentiable in $\mu$ such that \eqref{LL1} holds.

It remains to verify the $L$-differentiability.  By the chain rule and \eqref{LL2}, we obtain
\beg{align*} &F(P_{s,t}^* \mu\circ({\rm Id} +\phi)^{-1}) - F(P_{s,t}^*\mu)- D_\phi^L F(P_{s,t}^*\cdot)(\mu)
= \int_0^1\ff{\d}{\d\vv} F(\L_{X_{s,t}^\vv})\d\vv - D_\phi^L F(P_{s,t}^*\cdot)(\mu) \\
&= \int_0^1 \big\{\E\<(D^L F)(P_{s,t}^*\mu\circ({\rm Id}+\vv\phi)^{-1}) (X_{s,t}^\vv), v_{s,t}^{\phi,\vv}\>-  \E\<(D^L F)(P_{s,t}^*\mu)(X_{s,t}^\mu), v_{s,t}^{\phi,0}\>\big\}\d\vv.\end{align*}
Combining this with $F\in C^1(\scr P_2)$ with bounded $\|D^LF\|$,   the continuity of $v_{s,t}^{\phi,\vv}$ in $\vv$, \eqref{ES9}, and that $X_{s,t}^\vv\to X_{s,t}^\mu$ when $\|\phi\|_{L^2(\mu)}\to 0$, by the dominated theorem we prove
$$\lim_{\|\phi\|_{L^2(\mu)} \downarrow 0} \ff{|F(P_{s,t}^* \mu\circ({\rm Id} +\phi)^{-1}) - F(P_{s,t}^*\mu)- D_\phi^L F(P_{s,t}^*\cdot)(\mu)|}{\|\phi\|_{L^2(\mu)}}=0,$$
thus, $F(P_{s,t}^*\mu)$ is $L$-differentiable in $\mu$.

\end{proof}

\section{Proof of Theorem \ref{T1} }
According to Lemma \ref{L1}, \eqref{MA0} and \eqref{ES8}, to estimate $\|D^L P_{s,t}f(\mu)\|$, it suffices to handle
the   derivative of $ p_{s,t}^\mu$ and $H_{s,r,t}^{\mu,m}$ in $\mu$.  To this end, for fixed $T>0$, we introduce the Gaussian  heat kernel
 \beq\label{ES10} h_{T}(s,y)=\ff{\exp[-\ff{|y|^2}{8sK_T} ]}{(8\pi s K_T)^{\ff d 2}},\ \ y \in \R^d, s>0,\end{equation} which satisfies the Chapman-Kolmogorov equation
\beq\label{CK} \int_{\R^d} h_T(s_1, y-z) h_T(s_2, z)\d z= h_T(s_1+s_2, y),\ \ s_1,s_2>0, y\in \R^d.\end{equation}
By {\bf(H)}, there exists a constant $K_1(T)$, which increases in $T$, such that
 \begin{align*}  p_{s,r,t}^{\mu,z}(y,z)
&\le K_1(T) h_{T}(t-r, y-z) \e^{-\ff{|y-z|^2}{8(t-r)K_T}},\ \ y,z\in\R^d, 0\le s\le r<t\le T,\mu\in \scr P_2.\end{align*}
Consequently, there exists a constant $K_2(T)$, which increases in $T$, such that
\beq\label{G2}\beg{split} &p_{s,r,t}^{\mu,z}(y,z) \Big(1+ \ff{|y-z|^2}{t-r} +\ff{|y-z|}{(t-r)^{\ff 1 2}}\Big)\\
& \le K_2(T) h_{T}(t-r, y-z),\ \ y,z\in\R^d, 0\le s\le r<t\le T,\mu\in \scr P_2.\end{split}\end{equation}

\beg{lem}\label{L3} Assume {\bf (H)}.  There exists a constant  $\bar{K}_T>0$ which increases in $T>0$, such that  for any $ 0\le s\le r<t\le T, y,z\in\R^d$ and $m\ge 1$,
$p_{s,r,t}^{\mu,z}(y,z)$ and $H_{s,r,t}^{\mu,m} $ are $L$-differentiable in $\mu\in \scr P_2$ satisfying
\beg{align} & \label{Ta1} \|D^L\{ p_{s,r,t}^{\cdot,z}(y,z) \}(\mu) \|\le \bar{K}_T  h_T (t-r, y-z),\\
&\label{Ta2}  | H_{s,r,t}^{\mu,m}(y,z) |  \le \ff{\bar{K}_T^m (t-r)^{\ff m 2 -1} }{\GG(\ff m 2)}  h_T  (t-r, y-z),\ \ m\ge 1,\\
& \label{Ta3} \|D^L \{H_{s,r,t}^{\cdot,m}(y,z)\} (\mu) \|\le   \ff{m \bar{K}_T^m (t-r)^{\ff m 2 -1} }{\GG(\ff m 2)}  h_T (t-r, y-z),\ \ m\ge 1.\end{align}

\end{lem}
\beg{proof}   By {\bf (H)}, we have $|m_{s,r,t}^{\mu,z}|\le (t-r)K_T$, so that   \eqref{G2}  yields
\beq\label{ES14}\begin{split} & p_{s,r,t}^{\mu,z} (y,z) \Big(1+\ff{|y-z-m_{s,r,t}^{\mu,z}|^2} {t-r}+ \ff{|y-z-m_{s,r,t}^{\mu,z}|}{(t-r)^{\ff 1 2}} \Big) \le C_1(T)   h_T(t-r, y-z)  \end{split}\end{equation}
for some constant $C_1(T)>0$ increasing in $T$, and all $0\le s\le r<t\le T, \mu \in \scr P_2$ and $y,z\in \R^d$.
 Combining this with {\bf (H)},  \eqref{ES6}, \eqref{ES14} and applying Lemma \ref{L2}, we   prove the $L$-differentiability of $p_{s,r,t}^{\mu,z}(y,z)$ in $\mu\in\scr P_2$ and  the estimate \eqref{Ta1}.

Next,  by {\bf (H)}, \eqref{ES6}, \eqref{ES7} and \eqref{ES14}, we find   constants $C_2(T), C_3(T)>0$ increasing in $T>0$ such that for any  $0\le s\le r<t\le T, \mu \in \scr P_2$ and $y,z\in \R^d$,
\beq\label{A1} \beg{split} |H_{s,r.t}^\mu(y,z)|&\le C_2(T)p_{s,r,t}^{\mu,z} (y,z) |y-z|\Big(\ff {1}{t-r}+\ff{|y-z-m_{s,r,t}^{\mu,z}|^2}{ (t-r)^2}+ \ff{|y-z-m_{s,r,t}^{\mu,z}|}{ t-r} \Big) \\
&\le  C_3(T) (t-r)^{-\ff 1 2}  h_T(t-r, y-z).   \end{split}\end{equation}
Assume that for some $k\ge 1$ we have
$$ | H_{s,r,t}^{\mu,k}(y,z) |   \le C_3(T)^k    (t-r)^{\ff k 2 -1} \Big(\prod_{i=1}^{k-1} \bb\Big(\ff i 2, \ff 1 2\Big) \Big) h_T  (t-r, y-z).$$
Combining this with \eqref{HU}, \eqref{CK},  and \eqref{A1},  we derive
\beg{align*} & | H_{s,r,t}^{\mu,k+1}(y,z) | \le  \int_r^t\d u \int_{\R^d}  |H_{s,u,t}^{\mu, k}(z',z) H_{s,r,u}^{\mu}(y,z') | \d z'\\
&\le  C_3(T)^{k+1}  h_T(t-r, y-z)   \Big(\prod_{i=1}^{k-1} \bb\Big(\ff i 2, \ff 1 2\Big) \Big)  \int_r^t (t-u)^{\ff k 2-1} (u-r)^{-\ff 1 2}\d u\\
&=  C_3(T)^{k+1} (t-r)^{\ff{k+1}{2}-1}  h_T(t-r, y-z)   \Big(\prod_{i=1}^{k} \bb\Big(\ff i 2, \ff 1 2\Big) \Big).\end{align*}
In conclusion, for any $m\ge 1$, we have
$$| H_{s,r,t}^{\mu,m}(y,z) |  \le C_3(T)^m  (t-r)^{\ff m 2 -1}   \Big(\prod_{i=1}^{m-1} \bb\Big(\ff i 2, \ff 1 2\Big) \Big) h_T  (t-r, y-z),$$
which implies \eqref{Ta2}  for $\bar{K}_T= C_3(T)\GG(\ff 1 2)$,  since
\beq\label{AU} \prod_{i=1}^{m-1} \bb\Big(\ff i 2, \ff 1 2\Big)= \prod_{i=1}^{m-1} \ff{\GG(\ff i 2)\GG(\ff 1 2)}{\GG(\ff {i+1}2)} = \ff{\GG(\ff 1 2)^m}{\GG(\ff m 2)}.\end{equation}

Finally, by {\bf (H)}, \eqref{ES6'}, \eqref{ES7}, Lemma \ref{L2} and \eqref{A1}, we see that $H_{s,r,t}^{\mu,m}$ is $L$-differentiable in $\mu$, and there exist  constants  $C_4(T), C_5(T)\ge C_3(T)$ increasing in $T>0$ such that
\beq\label{A2} \beg{split}   &\| D^L\{H_{s,r,t}^\cdot (y,z)\}(\mu)\|\\
&\le C_4(T)p_{s,r,t}^{\mu,z} (y,z) |y-z|\Big(\ff {1}{t-r}+\ff{|y-z-m_{s,r,t}^{\mu,z}|^2}{ (t-r)^2}+ \ff{|y-z-m_{s,r,t}^{\mu,z}|}{ t-r} \Big) \\
&\le  C_5(T) (t-r)^{-\ff 1 2}  h_T(t-r, y-z).   \end{split}\end{equation}
Assume that for some $k\ge 1$ we have
$$\| D^L \{H_{s,r,t}^{\cdot,k}(y,z)\}(\mu)\|   \le k C_5(T)^k    (t-r)^{\ff k 2 -1} \Big(\prod_{i=1}^{k-1} \bb\Big(\ff i 2, \ff 1 2\Big) \Big) h_T  (t-r, y-z).$$
Combining this with \eqref{HU}, \eqref{CK},  and \eqref{A2},   we derive
\beg{align*} & \| D^L \{H_{s,r,t}^{\cdot,k+1}(y,z) \}(\mu)\|  \\
& \le   \int_r^t\d u \int_{\R^d} \Big\{\|D^L\{H_{s,u,t}^{\cdot, k}(z',z)\}(\mu)\|\cdot | H_{s,r,u}^{\mu}(y,z') |  \\
&\qquad\qquad  \qquad \quad+ |H_{s,u,t}^{\mu, k}(z',z)  | \cdot\| D^L\{H_{s,r,u}^{\cdot}(y,z')\}(\mu)\|\Big\}\d z'\\
&\le  (k+1) C_5(T)^{k+1}  h_T(t-r, y-z)   \Big(\prod_{i=1}^{k-1} \bb\Big(\ff i 2, \ff 1 2\Big) \Big)  \int_r^t (t-u)^{\ff k 2-1} (u-r)^{-\ff 1 2}\d u\\
&=  (k+1)C_5(T)^{k+1} (t-r)^{\ff{k+1}2-1}  h_T(t-r, y-z)   \Big(\prod_{i=1}^{k} \bb\Big(\ff i 2, \ff 1 2\Big) \Big).\end{align*}
This together with \eqref{AU} implies   \eqref{Ta3} for $\bar{K}_T=C_5(T)\GG(\ff 1 2). $
\end{proof}

We are now ready to prove the main result.

\beg{proof} [Proof of Theorem \ref{T1}] By Lemma \ref{L3} with \eqref{ES8} and \eqref{CK},   $p_{s,t}^\mu(x,z)$ is
$L$-differentiable in $\mu$ for $t>s$, and
there exists a constant $\dd_T>0$ increasing in $T>0$ such that
  \beg{equation}\begin{split} \label{Dp} & \|D^L \{p_{s,t}^\cdot(x,z)\}(\mu)\|\le \bar{K}_T h_T(t-s, x-z)\\&\qquad + \sum_{m=1}^\infty \ff{(m+1) \bar{K}_T^{m+1}}{\GG(\ff m 2)} \int_s^t (t-s)^{\ff m 2 -1} \d r \int_{\R^d} h_T(t-r, y-z) h_T(r-s, x-y)\d y\\  &\le \dd_T h_T(t-s, x-z). \ \ \  0\le s<t\le T, x,z\in \R^d, \mu\in \scr P_2.
  \end{split}\end{equation}
 This and \eqref{ES1} imply that $P_{s,t}^\mu f(x)$ is $L$-differentiable in $\mu$  such that
 $$\|D^L\{P_{s,t}^\cdot f(x)\}(\mu)\|\le \|f\|_\infty\int_{\R^d} \|D^L \{p_{s,t}^\cdot (x,z)\}(\mu)\|\d z\le \dd_T\|f\|_\infty$$
 holds for all $0\le s<t\le T, f\in \B_b(\R^d)$ and $\mu\in \scr P_2$. Combining this with Lemma \ref{L1}, we prove that $P_{s,t}f(\mu)$ is $L$-differentiable in $\mu$ and \eqref{EST1} holds for some increasing $C: [0,\infty)\to (0,\infty).$
According to the proof of \cite[Corollary 2.2(2)]{RW}, we can show that \eqref{EST1} implies \eqref{EST2}. We include below  a simple  proof   for completeness.

Since $C_b^1(\R^d)$ is dense in $L^1(P_{s,t}^*\mu+P_{s,t}^*\nu)$, \eqref{EST2} is equivalent to
 \beq\label{EST2'} |P_{s,t} f(\mu)- P_{s,t} f(\nu)|\le \ff{C_t\|f\|_\infty}{\sqrt{t-s}}\W_2(\mu,\nu),\ \ t>s, f\in C_b^1(\R^d), \mu,\nu\in \scr P_2.\end{equation}
Let $f\in C_b^1(\R^d)$ be fixed. We first prove this inequality  for     $\mu,\nu$ with compact supports.
 Let $\xi,\eta$ be two bounded random variables such that $\L_\xi=\mu, \L_\eta=\nu$ and
 $$\E |\xi-\eta|^2= \W_2(\mu,\nu)^2.$$ By  Proposition 3.1 in \cite{RW} and \eqref{EST1}, we obtain
\beg{align*}  &|P_{s,t}f(\mu)-P_{s,t}f(\nu)|= \bigg|\int_0^1 \ff{\d }{\d r} P_{s,t} f(\L_{r \xi+(1-r)\eta})  \d r\bigg| \\
&\le \int_0^1 \big|\E\<D^L P_{s,t}f(\L_{r \xi+(1-r)\eta})  (r \xi+(1-r)\eta), \xi-\eta\> \big|\d r \le \ff{C_t \|f\|_\infty}{\ss{t-s}} \W_2(\mu,\nu).\end{align*}
So, \eqref{EST2'} holds.

Next, for any $\mu,\nu\in \scr P_2$,  we choose   $\{\mu_n,\nu_n\}_{n\ge 1}\subset \scr P_2$ with compact supports such that
$$\lim_{n\to\infty}\big\{\W_2(\mu,\mu_n)+\W_2(\nu,\nu_n)\big\}=0.$$
Then  by the last step,
\beq\label{LST2} |P_{s,t}f(\mu_n)-P_{s,t}f(\nu_n)|\le \ff{C_t \|f\|_\infty}{\ss{t-s}} \W_2(\mu_n,\nu_n),\ \ n\ge 1.\end{equation} If $P_{s,t}f(\gg)$ is continuous in $\gg\in \scr P_2$, then by letting $n\to\infty$ we obtain the desired estimate \eqref{EST2'}.  To prove the continuity,
for any $\gg_1,\gg_2\in \scr P_2$,  let $\xi_1,\xi_2$ be $\F_0$-measurable random variables  such that $\L_{\xi_i}=\gg_i, i=1,2,$ and
$$\W_2(\gg_1,\gg_2)^2 =\E|\xi_1-\xi_2|^2.$$  For any   $\vv\in [0,1]$, let $X_{s,t}^\vv$ solve \eqref{E11} with initial value $X_{s,s}^\vv:= \vv\xi_1+(1-\vv)\xi_2$.
By \cite[Proposition 3.2 and (4.2)]{RW},
$$\nn_{\xi_1-\xi_2} X_{s,t}^\vv:=\ff{\d}{\d\vv} X_{s,t}^\vv$$ exists in $L^2(\P)$ with
$$\E|\nn_{\xi_1-\xi_2} X_{s,t}^\vv|^2\le c(t)\E|\xi_1-\xi_2|^2= c(t) \W_2(\gg_1,\gg_2)^2$$ for some constant $c(t)>0$.  Then
\beg{align*} &|P_{s,t}f(\gg_1)-P_{s,t}f(\gg_2)| =|\E f(X_{s,t}^1)-\E f(X_{s,t}^0)|  = \bigg|\int_0^1 \ff{\d }{\d \vv} \E f(X_{s,t}^\vv)    \d \vv \bigg| \\
&\le \int_0^1 \big|\E\<\nn f ( X_{s,t}^\vv), \nn_{\xi_1-\xi_2} X_{s,t}^\vv\> \big|\d s \le \ss{ c(t) }\|\nn f\|_\infty  \W_2(\gg_1,\gg_2).\end{align*}
Therefore, $P_{s,t}f(\gg)$ is continuous in $\gg\in\scr P_2$ and the proof is then finished.  \end{proof}

\end{document}